\def\draft{n}
\documentclass[headings]{amsart}
\usepackage{fullpage,amssymb,epsfig,amsbsy,amsmath,
pb-diagram}

\theoremstyle{plain}

\newtheorem{theorem}{Theorem}
\newtheorem{proposition}{Proposition}[section]

\theoremstyle{definition}
\newtheorem{definition}[proposition]{Definition}

\theoremstyle{remark}
\newtheorem{example}[proposition]{Example}

\newtheorem{remark}[proposition]{Remark}

\def\printname#1{
    \if\draft y
        \smash{\makebox[0pt]{\hspace{-0.5in}
            \raisebox{8pt}{\tt\tiny #1}}}
    \fi
}

\newcommand{\psdraw}[2]
         {\begin{array}{c} \hspace{-1.3mm}
    \raisebox{-4pt}{\epsfig{figure=draws/#1.eps,width=#2}}
    \hspace{-1.9mm}\end{array}}

\newlength{\standardunitlength}
\catcode`\@=11
\long\def\@makecaption#1#2{%
     \vskip 10pt

\setbox\@tempboxa\hbox{
       \small\sf{\bfcaptionfont #1. }\ignorespaces #2}%
     \ifdim \wd\@tempboxa >\captionwidth {%
         \rightskip=\@captionmargin\leftskip=\@captionmargin
         \unhbox\@tempboxa\par}%
       \else
         \hbox to\hsize{\hfil\box\@tempboxa\hfil}%
     \fi}
\font\bfcaptionfont=cmssbx10 scaled \magstephalf
\newdimen\@captionmargin\@captionmargin=2\parindent
\newdimen\captionwidth\captionwidth=\hsize
\catcode`\@=12

\newcommand{\tr}{\operatorname{tr}}

\def\lbl#1{\label{#1}\printname{#1}}

\def\BQ{\mathbb Q}

\def\A{\mathcal A}

\def\M{\mathcal M}

\def\H{\mathcal H}
\def\be{\begin{equation}}
\def\ee{\end{equation}}

\def\La{\Lambda}
\def\l{\lambda}

\def\S{\Sigma}

\def\s{\sigma}

\def\ga{\gamma}

\def\la{\langle}
\def\ra{\rangle}

\def\d{\delta}

\def\s{\sigma}

\def\sub{\subset}

\def\pt{\partial}

\def\Sym{\mathrm{Sym}}

\def\longto{\longrightarrow}

\def\hb{\hbar}

\def\bl{\boldsymbol{\lambda}}                   

\def\bs#1{\boldsymbol{#1}}

\def\fg{\mathfrak{g}}
\def\ZZ{\mathcal Z}
\def\FF{\mathcal F}
\def\ZU{Z^{\mathrm{U}}}

\def\FU{F^{\mathrm{U}}}
\def\glN{\mathfrak{gl}_N}
\def\deg{\mathrm{deg}}

\def\Stringp{\mathrm{Str}^{\mathrm{p}}}

\newdimen\tableauside\tableauside=1.0ex
\newdimen\tableaurule\tableaurule=0.4pt
\newdimen\tableaustep
\def\phantomhrule#1{\hbox{\vbox to0pt{\hrule height\tableaurule
width#1\vss}}}
\def\phantomvrule#1{\vbox{\hbox to0pt{\vrule width\tableaurule
height#1\hss}}}
\def\sqr{\vbox{%
  \phantomhrule\tableaustep

\hbox{\phantomvrule\tableaustep\kern\tableaustep\phantomvrule\tableaustep}%
  \hbox{\vbox{\phantomhrule\tableauside}\kern-\tableaurule}}}
\def\squares#1{\hbox{\count0=#1\noindent\loop\sqr
  \advance\count0 by-1 \ifnum\count0>0\repeat}}
\def\tableau#1{\vcenter{\offinterlineskip
  \tableaustep=\tableauside\advance\tableaustep by-\tableaurule
  \kern\normallineskip\hbox
    {\kern\normallineskip\vbox
      {\gettableau#1 0 }%
     \kern\normallineskip\kern\tableaurule}%
  \kern\normallineskip\kern\tableaurule}}
\def\gettableau#1 {\ifnum#1=0\let\next=\null\else
  \squares{#1}\let\next=\gettableau\fi\next}

\begin{document}


\title[On Chern-Simons Matrix Models]{
On Chern-Simons Matrix Models}

\author{Stavros Garoufalidis}
\address{School of Mathematics \\
         Georgia Institute of Technology \\
         Atlanta, GA 30332-0160, USA \\
         {\tt http://www.math.gatech} \newline {\tt .edu/$\sim$stavros } }
\email{stavros@math.gatech.edu}
\author{Marcos Mari\~no}
\address{Department of Physics, Theory Division \\
        CERN \\
         Geneva 23, CH-1211 Switzerland and Departamento de Matem\'atica,
IST, Lisboa, Portugal}
\email{marcos@mail.cern.ch}

\thanks{The first author was supported in part by NSF. \\
\newline
1991 {\em Mathematics Classification.} Primary 57N10. Secondary 57M25.
}

\date{March 25, 2005 \hspace{0.5cm} First edition: March 25, 2005.}

\begin{abstract}
The contribution of reducible connections to the $U(N)$ Chern-Simons
invariant of a Seifert manifold $M$ can be
expressed in some cases in terms of matrix integrals. We show that the
$U(N)$ evaluation of the LMO invariant
of any rational homology sphere admits a matrix model representation which
agrees with the Chern-Simons matrix integral for Seifert spheres and the trivial connection.
\end{abstract}

\maketitle

\tableofcontents


\section{Introduction}
\lbl{sec.intro}
Chern-Simons invariants of links and three-manifolds have been a rich arena for the interactions of
mathematics and physics in the last years. More recently, there has been a growing connection between
Chern-Simons invariants and topological string theory/Gromov-Witten theory. This has motivated various
developments and results. One of these developments has been the representation of Chern-Simons
invariants of three-manifolds in terms of matrix integrals over a Lie algebra \cite{M,AKMV}. This representation
has its origin in the work of Rozansky on the trivial connection contribution to the Chern-Simons invariant \cite{roz}, and
on the results of Lawrence and Rozansky on the $SU(2)$ Chern-Simons invariant of Seifert homology spheres \cite{LR}.

In this short note we clarify these results in the light of the LMO invariant \cite{LMO} and
its Aarhus integral representation \cite{aarhus}. After reviewing in section 2 the connection between
Chern--Simons theory and matrix integrals, we show in
section 3 that the LMO invariant of a rational homology sphere
$M$, evaluated for the $U(N)$ weight system,
can be always represented as a matrix integral. If the manifold $M$ is obtained through surgery on a link ${\mathcal L}$ in ${\bf S}^3$, the
matrix model `potential' involved in the matrix integral is related to the Kontsevich integral of ${\mathcal L}$. Since the LMO invariant
is conjectured to capture the trivial connection contribution to the Chern-Simons invariant of $M$, this suggests that this contribution
always has a matrix integral representation, as shown in \cite{M}.

\subsection{Acknowledgement}
M.M. would like to thank D. Bar-Natan and N. Wyllard for discussions.

\section{The Witten-Reshetikhin-Turaev invariant as a matrix integral}
\lbl{sec.fti}
\subsection{The Witten-Reshetikhin-Turaev invariant}

The Witten-Reshetikhin-Turaev (WRT) invariant of a three-manifold $M$ was
originally defined by Witten in \cite{cs} as the partition function
of a certain topological quantum field theory on $M$, the so-called
Chern-Simons theory. Fortunately, the invariant can be defined
in a purely combinatorial way, as we now describe.

The WRT invariant depends on a choice of gauge group $G$ and of an integer $k$ related
to the level of the affine Lie algebra based on $G$. We will use the following notations:
$r$ denotes
the rank of $G$, and $d$ its dimension.
$y$ denotes the dual Coxeter number.
The fundamental weights will be denoted by
$\lambda_i$, and the simple roots by $\alpha_i$, with $i=1, \cdots, r$.
$\rho$ denotes as usual the Weyl vector,
given by the sum of the fundamental weights.
The weight and root lattices of $G$ are
denoted by $\Lambda_{\rm w}$ and $\Lambda_{\rm r}$, respectively.
Finally, we put $l=k +y$.

The WRT invariant can be defined in terms of a surgery presentation of $M$.
By theorem due to Lickorish,
any three-manifold
$M$ can be obtained by surgery on a link ${\mathcal L}$ in ${\bf S}^3$.
Let us denote by
${\mathcal K}_i$, $i=1, \cdots, L$, the components of ${\mathcal L}$.
The surgery operation means that
around each of the knots ${\mathcal K}_i$ we take a tubular
neighborhood ${\rm Tub}({\mathcal K}_i)$ that we remove from ${\bf S}^3$. This
tubular neighborhood is a solid torus with a contractible cycle $\alpha_i$ and
a noncontractible cycle $\beta_i$. We then glue the solid torus back after
performing an
${\rm SL}(2, {\bf Z})$ transformation given by the matrix
\begin{equation}
\label{sltwo}
U^{(p_i,q_i)}=\left( \begin{array} {cc} p_i & r_i \\
q_i & s_i\end{array} \right).
\end{equation}
This means that the cycles $p_i \alpha_i + q_i \beta_i$ and
$r_i \alpha_i + s_i \beta_i$ on the boundary of the complement of ${\mathcal
K}_i$ are identified with the cycles $\alpha_i$, $\beta_i$ in ${\rm
Tub}({\mathcal K}_i)$.

To define the WRT invariant we use the representation of ${\rm SL}((2, {\bf Z})$
in the space of integrable representations of the affine Lie algebra associated
to $G$. A representation given by a highest weight $\Lambda$ is integrable if
the weight $\rho + \Lambda$ is in the
fundamental chamber ${\mathcal F}_l$ The fundamental chamber
is given by $\Lambda_{\rm w}/l
\Lambda_{\rm r}$ modded out by the action of the Weyl group.
In the following,
the basis of integrable representations will be labeled by the
weights in ${\mathcal F}_l$. In the case of simply-laced gauge groups,
the ${\rm Sl}(2, {\bf Z})$ transformation given by $U^{(p,q)}$ has
the following matrix elements in the above basis \cite{roz,ht}:
\begin{eqnarray}
{\mathcal U}^{(p, q)}_{\alpha \beta} &=&
{[i \, {\rm sign}(q)]^{|\Delta_+|} \over (l |q|)^{r/2}}
\exp \Bigl[ -{ i d \pi \over 12}  \Phi (U^{(p,q)})\Bigr]
\Biggl( { {\rm Vol}\, \Lambda_{\rm w} \over {\rm Vol}\,  \Lambda_{\rm r}}
\Biggr)^{1 \over 2} \nonumber\\
& \cdot& \sum_{n \in \Lambda_{\rm r}/q \Lambda_{\rm r}}\sum_{w \in {\mathcal W}} \epsilon
(w) \exp \Bigl\{ {i \pi \over l q} ( p \alpha^2 - 2\alpha (l n + w(\beta))
+ s(ln + w(\beta))^2 \Bigr\}.
\end{eqnarray}
In this equation, $|\Delta_+|$ denotes the number of positive roots of $G$,
and the second sum is over the Weyl group ${\mathcal W}$ of $G$.
 $\Phi (U^{(p,q)})$ is the Rademacher function:
\begin{equation}
\label{rade}
\Phi\left[ \begin{array}{cc} p & r \\ q&s\end{array} \right]=
{p + s \over q} - 12 s(p,q),
\end{equation}
where $s(p,q)$ is the Dedekind sum
\begin{equation}
s(p,q)={1 \over 4q} \sum_{n=1}^{q-1} \cot \Bigl( {\pi n \over q}\Bigr)
\cot \Bigl( {\pi n p\over q}\Bigr).
\end{equation}
In terms of the above data, the WRT invariant
of $M$ is given by:
\begin{equation}
\label{partf}
Z(M, l)= {\rm e}^{i\phi_{\rm fr}} \sum_{\alpha_1, \cdots, \alpha_L
 \in {\mathcal F}_l}Z_{\alpha_1, \cdots, \alpha_L}({\mathcal L})\,  {\mathcal U}^{(p_1,
 q_1)}_{\alpha_1 \rho} \cdots {\mathcal U}^{(p_L,
 q_L)}_{\alpha_L \rho}.
\end{equation}
In this equation, $Z_{\alpha_1, \cdots, \alpha_L}({\mathcal L})$ is the
quantum group invariant of the link ${\mathcal L}$ with representation $\alpha_i-\rho$
attached to its $i$-th component (recall that the weights in ${\mathcal F}_l$
are of the form $\rho + \Lambda$). The phase factor ${\rm e}^{i\phi_{\rm
fr}}$ is a framing correction that guarantees that the resulting invariant
is in the canonical framing for the three-manifold $M$. Its explicit
expression is
\begin{equation}
\label{framcorr}
\phi_{\rm fr}= {\pi k d \over 12 l }\biggl(\sum_{i=1}^L \Phi(U^{(p_i,
q_i)}) - 3 \sigma \, ({\mathcal L})\biggr),
\end{equation}
where $\sigma({\mathcal L})$ is the signature of the linking matrix of ${\mathcal
L}$. One can show that the above definition of WRT is invariant under Kirby moves,
therefore it defines a topological invariant of the three-manifold $M$.

The WRT invariant was originally defined by Witten in quantum-field theoretic terms,
as the partition function of Chern-Simons theory on the three-manifold $M$. The action
of Chern-Simons theory is given by
\begin{equation}
\label{csaction}
S_{\rm CS}(A)={k \over 4 \pi} \int_M {\rm Tr} \Bigl( A \wedge dA + {2 \over 3} A\wedge
A
\wedge A \Bigr),
\end{equation}
where $A$ is a $G$-connection on $M$, and the WRT invariant is given by
\begin{equation}
\label{wrt}
Z_k(M)= \int {\mathcal D} A {\rm e}^{i S_{\rm CS}(A)}.
\end{equation}
The description given above in terms of combinatorial data can be
derived from Chern-Simons theory in the context of canonical quantization.

\subsection{Matrix integral representation of the WRT invariant}

The fact that the WRT invariant is the partition function of a quantum field theory
suggests that it can be evaluated semiclassically as a sum over critical points of the
action. In the case of the Chern-Simons
functional (\ref{csaction}), the critical points are flat connections on $M$.
Each term in the sum is in turn an asymptotic, perturbative
expansion around the flat connection in powers of
the coupling constant of the model. Moreover, it can be argued that the
terms in this perturbative expansion contain important topological information
about the three-manifold $M$. For example, the one-loop contribution involves the
analytic torsion of the flat connection \cite{cs}, while the two-loop contribution
around the trivial flat connection turns out to be equal to the Casson-Walker invariant of
$M$ \cite{roz}.

From the point of view of combinatorial definition of the WRT
invariant, the properties that emerge in the asymptotic expansion in powers of the
coupling constant are far from being obvious.
In the case of Seifert spaces, it was shown in \cite{LR} that the WRT invariant
for gauge group $SU(2)$ can be
written as a sum of contour integrals and residues which correspond to
the contributions associated to the different flat connections. Some of the results of
\cite{LR} were generalized in \cite{M} to general simply-laced groups, where it
was shown that the contribution of reducible flat connections can be written as a
matrix integral. In \cite{bw}, Beasley and Witten have presented a very elegant
derivation of this matrix integral representation in the case of the
trivial flat connection, by using non-abelian localization. A similar result has been
recently obtained in \cite{bt}. 

 Seifert homology spheres
 can be constructed by performing surgery on a
link ${\mathcal L}$ in ${\bf S}^3$
with $n+1$ components, consisting on $n$ parallel and unlinked unknots
together with a single unknot whose linking number with each of the other
$n$ unknots is one. The surgery data are $p_j/q_j$ for the unlinked
unknots, $j=1, \cdots, n$, and 0 on the final component. $p_j$ is coprime
 to $q_j$ for all $j=1, \cdots, n$, and the $p_j$'s are pairwise coprime.
After doing
surgery, one obtains the Seifert space
$M=X({p_1 \over q_1}, \cdots, {p_n
\over q_n})$. This is rational homology sphere whose
first homology group $H_1(M, {\bf Z})$ has order $|H|$, where
\begin{equation}
\label{orderh}
H=P \sum_{j=1}^n {q_j \over p_j}, \,\,\,\,\, {\rm and}
\,\,\,\,
P=\prod_{j=1}^n p_j.
\end{equation}
We will denote $e=H/P$. For $n=1,2$, Seifert homology spheres reduce to lens
spaces, and one has that $L(p,q)=X(q/p)$. For
$n=3$, we obtain the Brieskorn homology spheres $\Sigma (p_1, p_2, p_3)$
(in this case the manifold is independent of $q_1, q_2 ,q_3$).
By using the formulae for the WRT invariant presented above, one can write
the contribution of reducible flat connections to the
Chern-Simons partition function of $X({p_1 \over q_1}, \cdots, {p_n
\over q_n})$ as
\begin{eqnarray}
\label{fians}
& & {(-1)^{|\Delta_+|} \over | {\mathcal W}|\, (2 \pi i )^r}
\Biggl( { {\rm Vol}\, \Lambda_{\rm w} \over {\rm Vol}\,  \Lambda_{\rm r}}
\Biggr){ [{\rm sign}(P)]^{|\Delta_+|} \over |P|^{r/2}}
{\rm e}^{{ \pi i d \over 4} {\rm sign}(H/P) +{\pi i dy \over 12 l} \phi}
\nonumber\\
& \cdot & \sum_{t \in \Lambda_{\rm r}/H\Lambda_{\rm r}}
\int d\beta \,  {\rm e}^{ -\hbar e {\beta^2/2 } - l t\cdot \beta}{
\prod_{i=1}^n \prod_{\alpha>0} 2 \sinh {\beta \cdot \alpha \over
2 p_i } \over \prod_{\alpha >0}
\Bigl( 2 \sinh { \beta \cdot \alpha \over
2} \Bigr)^{n-2}}
\end{eqnarray}
In this equation, $\beta$ is an element of $\Lambda_w \otimes {\bf R}$, $\phi$ is given by
\be
\phi=e- 3 {\rm sign}\, (e)-12 \sum_{i=1}^n s(q_i,p_i).
\ee
and we have introduced
$$\hbar ={2\pi i \over k+y}.$$
The lattice $\Lambda_r/H \Lambda_r$ decomposes in different Weyl orbits, and each of these orbits
correspond to a different, reducible flat connection. The contribution of the trivial flat connection is obtained
by setting $t=0$ in (\ref{fians}).

\section{Matrix integrals and the LMO invariant}
\lbl{sec.lmo}

In this section, we will show that the LMO invariant of a rational homology sphere, evaluated
through the $U(N)$ weight system, can be always expressed as a matrix integral. This follows very simply by the
definition of the LMO invariant
in terms of formal Gaussian integration given in \cite{aarhus}, and the detailed structure of the $U(N)$
weight system. Since the LMO invariant is conjectured to capture the contribution of the trivial
connection to the WRT invariant,
our result indicates that this contribution is expected to have a representation in terms of matrix
integrals, as it happens
with the Seifert homology spheres. In particular, we will show that the result (\ref{fians}) for
Seifert spheres agrees with the $U(N)$ evaluation of the
LMO invariant calculated in \cite{bl}.

\subsection{A review of the Kontsevich integral}
\lbl{sub.kontsevich}

The physics origin of the Kontsevich integral of a link in ${\bf S}^3$ is
Chern-Simons perturbation theory along the trivial flat connection of the
backround 3-space. The Feynmann diagrams of the theory are trivalent graphs
with legs (so-called unitrivalent graphs). The legs are colored by the
components of the link, and the edges along the trivalent vertices are
equipped with a cyclic ordering. The graphs are considered modulo the
AS and IHX relations. The graphs can be multiplied (using the disjoint union)
and can be graded (where the degree of a graph is half the number of
vertices). Using formal linear combinations
(with coefficients in $\BQ$) of these graphs, we can define a completed
graded algebra $\A(\star_X)$
where $X$ is a set in 1-1 correspondence with the components of the link.

It turns out that the Konstevich integral of a link is a group-like element
of $\A(\star_X)$, thus we can define its logarithm
$$
\FF(S^3,L)=\log \ZZ(S^3,L)
$$
which lies in the completed vector space $\A^c(\star_X)$ generated by
{\em connected} unitrivalent graphs, modulo the AS and IHX relations.

There are two degrees of a diagram $D$ in $\A(\star_X)$:
\begin{itemize}
\item
The {\em Vassiliev degree}
$\deg_1(D)$, which equals to half the number of vertices.
\item
The {\em Euler degree} $\deg_2(D)$ which equals to $-\chi(D)$.
\end{itemize}

Notice that the Euler degree of a connected diagram is $\text{rk}H_1(D)-1
\geq -1$,
and that $\deg_1(D)=\deg_2(D)+ | \text{Legs(D)}|$.

\subsection{A review of  the LMO invariant}
\lbl{sub.LMOreview}

In this section we review the LMO integral. The physics origin of the
LMO invariant (and its cousin, the Aarhus integral) is Chern-Simons
perturbation theory along the trivial flat connection.

Consider a framed link $L$ of $r$ components in $S^3$, and let $M$ denote the 3-manifold
obtained by surgery on $L$.

For every nonnegative natural number $m$, there is an LMO integration map
$$
\int^{(m)} \, dX: \A(\star_X) \longto \A(\emptyset)
$$
which takes values in $\A_{\leq m}(\emptyset)$, the Vassiliev degree $\leq m$
part of the completed vector space $\A(\emptyset)$ of
trivalent graphs with vertex orientations modulo the AS and IHX relations.
If the integrand $a$ is group-like, then the sequence
$\{\int^{(m)} a \, dX\}_m$ can be assembled in a group-like element in
$\A(\emptyset)$, which we denote by $\int a\, dX$.

This allows us to define
$$
\ZZ_0(L)=\int \ZZ(S^3,L)
$$
as well as
$$
\ZZ(M)=\frac{\ZZ_0(L)}{\ZZ_0(S^3,U^+)^{\s^+(L)}
\ZZ_0(S^3,U^-)^{\s^-(L)}}
$$
where $U^{\pm}$ is the unknot with framing $\pm 1$ and $\s^{\pm}(L)$
is the number of positive (resp. negative) eigenvalues of the linking
matrix of $L$.
It turns out that $\ZZ(M)$ is depends only on $M$ and not on the framed link
$L$. Moreover, $\ZZ(M)$ is group-like, that is we can define its logarithm:
$$
\FF(M)=\log \ZZ(M)
$$
which takes values in $\A^c(\emptyset)$.

Here is a rough description of the LMO integration $\int^{(m)}$.
\begin{itemize}
\item
Consider an element $a \in \A(\star_X)$. Concentrate on $a_{2m}$, the
piece of $a$ that contains diagrams with exactly $2m$ legs of each color.
\item
Then, $\int^{(m)} a \, dX$ is the sum of all $((2m-1)!!)^{|X|}$ ways of
pairing up the legs of each color $X$.
\item
We consider the result in $\A_{\leq m}(\emptyset)$.
\item
If $a$ is group-like, then we can assemble the pieces $\int^{(m)} a \, dX$
for all $m$ into a group-like element in $\A(\emptyset)$.
\end{itemize}

\subsection{Weight systems}
\lbl{sub.weight}

For every semisimple Lie algebra $\fg$, we have a {\em weight system} map:
$$
W_{\fg}: \A^c(\star_X) \longto \hb S(\fg^{\oplus |X|})^{\fg}[[\hb]]
$$
Specifically, if $D$ is a unitrivalent graph with $2m$ vertices and $l$ legs
then $W_{\fg}(D) \in S^l(\fg^{\oplus |X|})^{\fg} \, \hb^l$. Notice that
$m=-\chi(D)+l$.


Let $W^U=W_{\glN}$, for arbitrary $N$. We now describe in detail the $W^U$ weight system. Following Bar-Natan, let us
introduce the vector space spanned by marked
surfaces.

\begin{definition}
\lbl{def.markedS}
An $X$-{\em marked surface} $(\S,\ga)$ is an oriented compact
topological surface $\S$ with nonempty boundary, together
with a choice $\ga$ of points (colored by $X$) on $\pt \S$.
Let $\M_X$ denote the completed vector space of formal $\BQ$-linear
combinations of $X$-marked surfaces.
\end{definition}

If $|X|=1$ and $\S$ is connected, $\ga$ gives rise to a partition
$(0^{\ga_0} 1^{\ga_1}
\dots)$, where $\ga_j$ is the number of boundary components of $\S$
with $j$ points.

Like unitrivalent graphs, marked surfaces have two degrees:
\begin{itemize}
\item
The {\em Vassiliev degree} $\deg_1(\S,\ga)$ of a marked surface is
$-\chi(\S)+|\ga|$, where $|\ga|=\sum_{j} j \ga_j$.
\item
The {\em Euler degree} $\deg_2(\S,\ga)$ of a marked surface is $-\chi(\S)$.
\end{itemize}

Thus, $|\ga|$ in the case of a marked surface plays the role of the number
of legs.

Marked surfaces can be multiplied (via the disjoint union) and graded
(by the Vassiliev degree). Let $\M^c_X$ be defined analogously.

There is a map:
$$
\Psi: \A(\star_X) \longto \M_X
$$
defined by:
$$
D \longto \sum_M (-1)^{s_M} (\S_{D,M},\ga_{D,M})
$$
where
\begin{itemize}
\item
the sum is over all possible markings $M$ of the trivalent vertices
of $D$ by $0$ or $1$,
\item
$s_M$ is the sum, over the set of trivalent vertices, of the values of $M$
\item
$\S_{D,M}$ denotes the $X$-marked surface
obtained by thickening the trivalent vertices of $D$ as follows:
\begin{equation}
\lbl{eq.thicken}
\psdraw{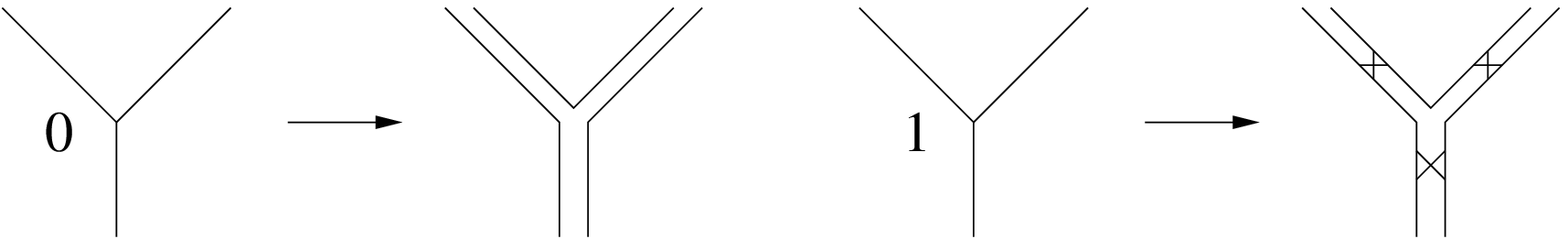}{2in}
\end{equation}
and thickening the edges of $D$, and connected up to a surface. The legs
of $D$ become the choice $\ga_{D,M}$ of points in $\S_{D,M}$.
It turns out that $\S_{D,M}$ is well-defined and oriented.
\end{itemize}

The above map preserves the Vassiliev and Euler degrees $\deg_1$ and $\deg_2$.

Moreover, there is a map:
$$
\Phi: \M^c_X \longto \La^{\otimes r}[N,\hb]
$$
where $\Lambda$ is the ring of symmetric polynomials. This map is
defined by:
$$
(\S,\ga) \longto
N^{\ga_0} \prod_{n=1}^\infty p_n^{\ga_n}
 \, \hb^{m} \in \La^{\otimes r}_l[N,\hb]
$$
where $m=\deg_1(\S,\ga)$, $l=\deg_1(\S)-\deg_2(\S)$, and $p_n$ is the power sum $\sum_{j}x_j^n$. We
remind that products of power sums provide a basis for the ring of symmetric polynomials in the variables $x_j$.
We will also use in the following the basis of $\Lambda$ given by Schur functions $s_{\lambda}$, which are labeled by a
partition $\lambda$. An $r$-uple of partitions will be denoted by ${\bf \lambda}=(\lambda_1, \cdots, \lambda_r)$. A
basis of $\Lambda^{\otimes r}$ is therefore provided by the products $s_{\bf \lambda}=s_{\lambda_1} \cdots s_{\lambda_r}$.


\begin{remark}
\lbl{rem.iso}
In fact, the map $\Phi$ is a vector space isomorphism, although
we will not use this.
\end{remark}

\begin{proposition}
\lbl{prop.UN}\cite{BN,CD}
We have a commutative diagram
$$
\begin{diagram}
\node{\A^{\text{gp}}(\star_X)}
\arrow[2]{e,t}{\Psi} \arrow{se,l}{W^U}
\node[2]{\exp(\M^c_X)}
\arrow{sw,r}{\Phi}   \\
\node[2]{\exp\left(\Stringp_r \right)}
\end{diagram}
$$
where we define
\begin{equation}
\lbl{eq.stringp}
\Stringp_r:=
\begin{cases}
\left\{\sum_{g=0}^\infty \sum_{\bl \neq \bs{0}}
a^g_{\bl} \, s_{\bl} \,\, \Big|
a^g_{\bl} \in \hb^{2g-2+|\bl|} \BQ[N,\hb] \right\}
\sub  \La^{\otimes r}[N,\hb] & \text{if} \quad r > 0 \\
\frac{1}{\hb^2} \BQ[N,\hb]  & \text{if} \quad r=0.
\end{cases}
\end{equation}
\end{proposition}

\begin{proof}
Let $e_{i,j}$ for $1 \leq i,j \leq N$ be a basis for
$\glN$. We have

\begin{eqnarray*}
\tr (e_{i,j} e_{k,l}) &=& \d_{i,l} \, \d_{j,k} \\
 \text{$[$} e_{i,j}, e_{k,l} \text{$]$} &=& \d_{j,k} \, e_{i,l}-\d_{l,i} \, e_{k,j}.
\end{eqnarray*}

A diagram consists of a number of $Y$ graphs some of whose half-edges are
glued in pairs. The weight system colors each half-edge by an element
of $\glN$. Graphically, the above equations become:
$$
\psdraw{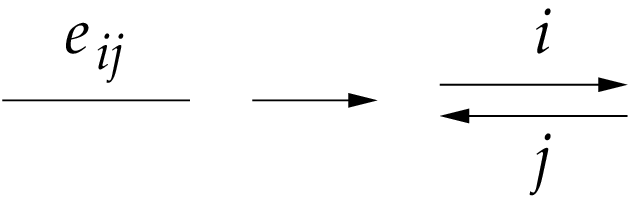}{1in}
$$
and Equation \eqref{eq.thicken}.

This computes the corresponding element in $S(\glN)^{\glN}$.

Let $E=(e_{ij})$ denote the $N$ by $N$ matrix with noncommutative
entries $e_{ij}$.
Consider the $X$-marked ribbon graph of genus $0$:
$$
R_n:= \sum_{i_1,\dots,i_n=1}^N \psdraw{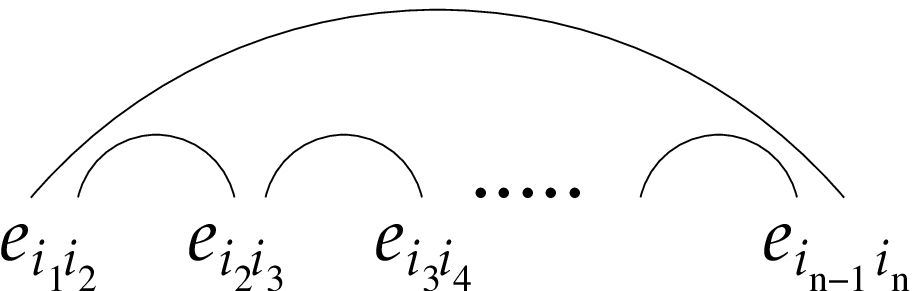}{2in}
$$
Then, it is easy to see that $W_{\glN}(R_n)=\tr(E^n)$.

We have $S(\glN)^{\glN} \cong S^(\mathfrak{h}_N)^{\Sym_N}$, where
$\mathfrak{h}_N$ is the Cartan subalgebra spanned by $x_i:=E_{ii}$ for
$i=1,\dots, N$. Under this isomorphism, $E$ maps to a diagonal matrix
$\text{diag}(x_1,\dots,x_N)$, thus $\tr(E^n)$ gets mapped to
$\sum_{i=1}^N x_i^n=p_n$.
\end{proof}

\begin{example}
\lbl{ex.wheel}
If $w_2:=\psdraw{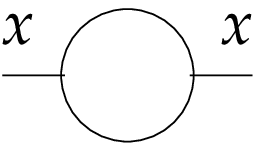}{0.3in}$ is the wheel with $2$ legs colored by
$X=\{x\}$, we have:
$$
w_2
 \stackrel\Psi\longto
 2 \left( \psdraw{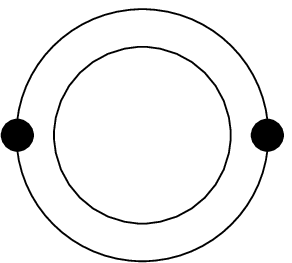}{0.5in} - \psdraw{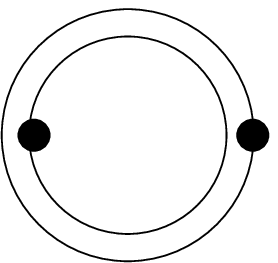}{0.5in} \right)
\stackrel\Phi\longto 2(N p_2-p_1^2) \hb^2
$$
For a polywheel $w_n$ with $n$ legs colored by $\{x\}$, we obtain
$$
(\Phi \, \Psi)(w_n)=\sum_{1\le i,j\le N} (x_i- x_j)^n=\hbar^n \sum_{s=0}^{2n} (-1)^s {2 n \choose s} p_s p_{2n-s},
$$
where we set $p_0=N$.
\end{example}

We will define
\begin{equation}
\lbl{eq.defineU}
\ZU:=W^U \, \ZZ \qquad
\FU:=W^U \, \FF .
\end{equation}
We may think of $\FU({\bf S^3},L)$ as the potential for a matrix model. As we will see, this
potential the 3-manifold invariant,
by integration on the full Lie algebra.





\subsection{Weight systems commute with LMO integration}
\lbl{sub.behavior}

\begin{definition}
\lbl{def.intX}
Let us define
$$
\int^{(m)}\, dX:  \M_X  \longto \M_{\emptyset}
$$
as follows:
$$
\int^{(m)} \prod_j p_j^{\ga_j}   \, dX=
\begin{cases}
\left\la
\sqcup_i (\sqcup^{\ga_i} R_i), \prod_{x \in X}
\frac{1}{m!} \left( \frac{\psdraw{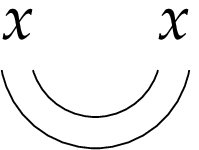}{0.3in}}{2} \right)
\right\ra_X & \text{if} \quad |\ga|=2m \\
0 & \text{if} \quad |\ga| \neq 2m
\end{cases}
$$
where for two ribbon surfaces $\S$ and $\S'$ with $X$-marked boundary,
$\la \S, \S' \ra_X$ is the sum over all pairings of the $X$-marked ends
of $\S$ with those of $\S'$ (if they match, otherwise zero).
\end{definition}

\begin{example}
\lbl{ex.pbasis}
For $m=1$, we have:
\begin{eqnarray*}
\int^{(1)} p_2 \, dX &=&
\left\la R_2, \frac{\psdraw{lowercupx}{0.3in}}{2} \right\ra=
\left\la \psdraw{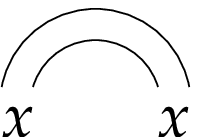}{0.3in}, \frac{\psdraw{lowercupx}{0.3in}}{2}
\right\ra=
\psdraw{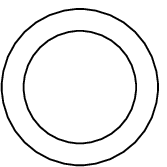}{0.3in} = N^2 \\
\int^{(1)} p_1^2 \, dX &=&
\left\la R_1 \sqcup R_1, \frac{\psdraw{lowercupx}{0.3in}}{2} \right\ra=
\left\la \psdraw{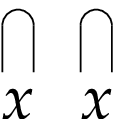}{0.15in}, \frac{\psdraw{lowercupx}{0.3in}}{2} \right\ra=
\psdraw{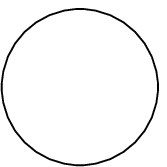}{0.3in} = N
\end{eqnarray*}
and
$$
\int^{(1)} w_2 \, dX= 2(N^3-N)
$$
\end{example}

\ref{def.intX} leads naturally to a map
\be
\label{intu}
\int^U: \exp \left(\Stringp_r\right) \rightarrow \left(\Stringp_0\right)
\ee
such that the following diagram commutes:
$$
\begin{diagram}
\node{\A^{\text{gp}}(\star_X)}
\arrow{e,t}{\Psi}\arrow{s,l}{\int}
\node{\exp\left(\M_X^c\right)}
\arrow{s,r}{\int} \arrow{e,t}{\Phi}
\node{\exp\left(\Stringp_r\right)}
\arrow{s,r}{\int^U}\\
\node{\A^{\text{gp}}(\emptyset)}
\arrow{e,t}{\Psi}
\node{\exp\left(\M_{\emptyset}^c\right)}
\arrow{e,t}{\Phi}
\node{\exp\left(\Stringp_0\right)}
\end{diagram}
$$


\subsection{The $U(N)$-version of the LMO invariant as a matrix model}
\lbl{sub.LMOmm}

The next theorem identifies the $\ZU(M)$ invariant of a closed 3-manifold
with a matrix model. Consider a framed link $L$ in $S^3$, and let $M$ denote
the closed 3-manifold resulting from Dehn surgery on $L$. Then, the image
$\ZU(S^3,L)$ of the Kontsevich integral of $L$ under the $\glN$ weight system
lies in $\exp(\Stringp_r)$, where $r=|L|$. We want to show that the integration
(\ref{intu}) induced by \ref{def.intX} agrees with Gaussian matrix model integration.

We first define Gaussian matrix model integration.

\begin{definition}
\lbl{def.intG} The map
$$
\int_{\H_N}: \exp \left(\Stringp_1\right) \rightarrow \left(\Stringp_0\right)
$$
is defined as follows:
$$
\int_{\H_N} \, p_{\l} ={1\over Z} \int dM \, \prod_{j} ({\rm tr}\, M^j)^{k_j}
 \exp\left(-\frac{1}{2} \tr(M^2)\right),
$$
where $\l=(1^{k_1}2^{k_2} \cdots)$, $M$ is a Hermitian, $N\times N$ matrix, and
$$
Z= \int dM \, \exp\left(-\frac{1}{2} \tr(M^2)\right).
$$
The measure $dM$ in the matrix integral is given by
$$
dM=\prod_{i=1}^N d M_{ii} \prod_{1\le i<j\le N} d({\rm Re}\, M_{ij}) d({\rm Im}\, M_{ij}).
$$
The above definition can be extended to $\exp \left(\Stringp_r\right)$ as follows
$$
\int_{\H_N^r} p_{\bs{\mu}} = \int_{\H_N}
p_{\mu_1} \times \dots \times
\int_{\H_N} p_{\mu_r}.
$$
\end{definition}

\begin{theorem}
\lbl{thm.LMOmm}
We have:
$$
\int^U = \int_{\H_N^r}
$$
As a result, the $\glN$ version of the LMO invariant is given by a matrix
model:
$$
\int_{\H_N^r} dX \, \ZU(S^3,L) = \ZU(M).
$$
\end{theorem}

This theorem follows immediately from the fact that \ref{def.intX} is simply the description
of Gaussian integration in terms of Wick contractions, see for example \cite{IZ, BIZ}.

\begin{example}
\lbl{ex.sbasis}
We have:
\begin{eqnarray*}
\int_{\H_N} s_2 \, dX &=& \frac{N(N+1)}{2} \\
\int_{\H_N} s_{1,1} \, dX &=& -\frac{N(N-1)}{2}
\end{eqnarray*}
and $p_2=s_2-s_{1,1}$ and $p_1^2=s_2+s_{1,1}$.
Comparing with Example \ref{ex.pbasis},
we confirm the above claim for partitions with $2$ boxes.
\end{example}

\begin{example}
The result of \cite{bl} expresses the LMO invariant of a Seifert sphere as
\be
Z =\exp\bigl( {\theta\over 48} (e_0 - 3 {\rm sign}(e_0) -\sum_i S({q_i\over p_i})) \bigr) \int^{(m)} \Omega^{2-n}_{x/e_0^{1/2}}\prod_{\ell=1}^n \Omega_{x/(e^{1/2}_0 p_\ell)}
\ee
where $\Omega_x$ is the element in  $\A(\star_X)$ (with $X=\{x\}$) introduced in \cite{wheels} and given by
$$
\Omega_x=\exp \sum_{m=1}^{\infty} b_{2m} w_{2m},
$$
where
$$
\sum_{m=0}^{\infty} b_{2m} x^{2m} ={1\over 2} \log {\sinh x/2 \over x/2}.
$$
One easily calculates
$$
(\Phi \, \Psi)(\Omega^{2-n}_{x}\prod_{\ell=1}^n \Omega_{x/p_\ell})=P^{|\Delta_+|} \Delta^{-2}(x) \prod_{i<j}
\Bigl(2 \sinh ({x_i-x_j\over 2}) \Bigr)^{2-n} \prod_{\ell=1}^n \prod_{i<j} \Bigl(2 \sinh
({x_i-x_j\over 2 p_{\ell}})\Bigr),
$$
where $\Delta(x)=\prod_{i<j}(x_i-x_j)$. On the other hand, it is well known that
Gaussian integration can be expressed in terms of
eigenvalues as (see for example \cite{BIZ})
$$
\int_{\H_N} \, p_{\l} ={ \int \prod_{i=1}^N dx_i  e^{-x_i^2/2} \Delta^2(x) \prod_j  p_j^{k_j}
 \over \int \prod_{i=1}^N dx_i  e^{-x_i^2/2}  \Delta^2(x)}.
 $$
After writing $\beta =\sum_i x_i e_i$, where $e_i$ is an orthonormal basis in $\Lambda_w$, we
find that $Z^{\rm U}$ is indeed given by the matrix integral (\ref{fians}), up to an overall factor
\end{example}

\ifx\undefined\bysame
    \newcommand{\bysame}{\leavevmode\hbox
to3em{\hrulefill}\,}
\fi

\end{document}